\documentclass[12pt]{elsarticle}

\usepackage{amsmath}
\usepackage{amsfonts}
\usepackage{theorem}
\usepackage{enumerate}
\usepackage{graphicx}
\usepackage{amssymb}
\usepackage{color}
\usepackage{lineno}
\usepackage{bm}

\usepackage{epsfig}
\usepackage{fancybox}
\usepackage{color}

\theoremheaderfont{\scshape}
\theorembodyfont{\upshape}

\theoremstyle{break} 
\theoremstyle{break} 

\newcommand {\QED}{\hfill$\blacksquare$}

\usepackage[normalem]{ulem}
\newcounter{corr}
\definecolor{violet}{rgb}{0.580,0.,0.827}
\newcommand{\corr}[3]{\typeout{Warning : a correction remains in page
\thepage}
				\stepcounter{corr}        
				{\color{red}\ifmmode\text{\,\sout{\ensuremath{#1}}\,}\else\sout{#1}\fi} 
       {\color{blue}#2} 
       {\color{violet} #3}} 

\begin{document}

\title{An Arbitrary-Order Moving-Mesh Finite Element Algorithm for One-Dimensional Implicit Moving Boundary Problems}

\author[1]{Matthew E Hubbard\corref{cor1}}
\ead{matthew.hubbard@nottingham.ac.uk}
\affiliation[1]{organization={University of Nottingham},addressline={School of Mathematical Sciences, University Park},city={Nottingham},postcode={NG7 2RD},country={UK}}

\author[2]{Thomas J Radley}
\ead{tjradleyresearch@gmail.com}
\affiliation[2]{organization={University of Montpellier},addressline={CNRS},city={Montpellier},country={France}}

\cortext[cor1]{Corresponding author}


\date{\today}



\begin{abstract}
 We present a one-dimensional high-order moving-mesh finite element method for moving
boundary problems where the boundary velocity depends implicitly on the solution
in the interior of the domain. The
method employs a conservative arbitrary Lagrangian-Eulerian approach to predict the evolution of the
approximate solution. It retains the order of accuracy of the underlying fixed-mesh
method by (i) computing the boundary velocities using a high-order approximation to
a distributed local conservation principle which generates a Lagrangian `flow' velocity
and (ii) assuming a continuous piecewise linear variation of the mesh velocity field in the
interior of the computational domain.

 Within each time-step the algorithm consists of two stages: the computation of a velocity field
with which to move both the domain boundary and the internal mesh, and the approximation
of the solution to the PDE on the updated mesh. Both internal and boundary velocities are
generated within the same framework, though it would be simple to replace the internal
mesh velocity field by applying a more traditional mesh movement strategy (for example, one
which seeks to equidistribute an error indicator at each time-step): the high-order accuracy
should be retained as long the discrete velocity which is used in the solution update
is assumed to vary linearly within mesh elements.
The proposed method is explicit in time, requires only the inversion of linear systems
of equations within each time-step, and is implemented fully in the physical domain,
so no mapping to a reference domain is employed. It attains arbitrary-order
accuracy in both space and time without the need to satisfy a discrete geometric conservation law.

 Numerical results are presented to demonstrate high-order approximation of both
boundary position and internal solution for three distinct systems which test the method in different ways: linear diffusion
with phase change across a moving internal interface; linear diffusion with absorption
generating a shrinking domain; nonlinear diffusion generating a growing domain with
finite boundary velocity.

\end{abstract}

\begin{keyword}
 moving meshes \sep finite elements \sep high-order accuracy \sep parabolic PDEs \sep arbitrary Lagrangian-Eulerian methods
\end{keyword}

\maketitle

\section{Introduction}
\label{sec:intro}

 Mesh movement can be applied to mesh-based algorithms for approximating
systems of partial differential equations (PDEs) to accommodate the movement of a boundary of
the computational domain or to reduce the computational time and memory required to generate
an approximation with a specified level of accuracy. When used for the latter purpose, moving
the mesh ($r$-adaptivity) provides an alternative to local adjustment of the approximation
space by refining or coarsening mesh elements ($h$-adaptivity) or by changing the degree of
the polynomial approximation ($p$-adaptivity). Mesh movement is the least commonly used of
these three techniques, but is well suited to fitting and following sharply defined moving
interfaces or internal features which can be most effectively resolved on anisotropic meshes
which remain aligned with these features. Furthermore, when the mesh can evolve without
significant deterioration in quality (or, in extremes, tangling), there is no need to change
the mesh connectivity through, for example, edge or face swapping. In contrast to $h$- and $p$-adaptivity,
there is no additional overhead incurred by keeping track of the degrees of freedom and, in parallel
computations, no need to adjust the partitioning.

 Ideally, it would be possible to select an optimally efficient combination of $h$-, $p$-
and $r$-adaptivity for a given simulation. In particular, $h$- and/or $p$-refinement would
be essential if the desired accuracy could not be obtained using the fixed number of
degrees of freedom given initially to the moving-mesh simulation.
 Research on $hp$-adaptivity is relatively mature, see for example \cite{CDGH17,Dem06,Dem08}, but the extension
of fixed-mesh algorithms to moving-mesh algorithms which demonstrate the same convergence
rates as the spatial mesh size $h$ and the time-step length $\Delta t$
has lagged behind for $p > 1$, even when the domain boundary is fixed.
A natural approach would be to consider a space-time formulation, \emph{e.g.}\ \cite{HR20},
but the mesh movement distorts the space-time prism elements that would normally be used,
producing non-standard element shapes. These are dealt with in \cite{HR20} by dividing them up
to give a $d+1$-dimensional mesh of simplices for each space-time slab, though a more flexible alternative might be to
apply more recent work on polytopic meshes (see, for example, \cite{CDGH17}). Such space-time methods
can be of arbitrary order in space and time, be unconditionally stable in time and
automatically satisfy the geometric conservation law (GCL) \cite{FGG01,TL79,TT61}: the
penalty for this is the requirement to solve robustly a very large (often nonlinear) system at
each time-step. This is particularly challenging when the boundary movement, and hence the mesh
at the new time level, depends on the solution.

 It is more common to use an arbitrary Lagrangian-Eulerian (ALE) framework as an alternative to
space-time methods. See, for example, \cite{DHPR04} for an outline of the basic principles of ALE
methods and a summary of key early developments.
 The ALE formulation allows the use of standard time-stepping schemes and, crucially for our work,
provides more flexibility when the mesh at the new time level depends on the approximate solution
of the PDE.
 In the late 1990s and early 2000s multiple time-stepping schemes were proposed which retained
second-order accuracy on moving meshes \cite{FN04,GGF03} and
satisfied a discrete geometric conservation law (DGCL) \cite{BG04,FGG01,LF96}.

 The first works (to the authors' knowledge) that demonstrate higher than second-order
accuracy in both space and time for the ALE formulation of a multidimensional moving-mesh
algorithm appeared in 2009 \cite{EGP09,PBP09}. The focus of \cite{EGP09} is on high-order
accuracy in time (and the role of the GCL in this), but the authors use P2-P1 Taylor-Hood finite elements
to eliminate the spatial error completely for a 2D incompressible Navier-Stokes test case in which
velocity varies quadratically in space and the mesh moves according to a smooth, predefined velocity
field within a fixed domain. This suggests third-order accuracy in space along with the fifth-order
accuracy demonstrated in time.
 The work in \cite{PBP09} concentrates on spatial accuracy (with similar emphasis on the role of the GCL),
using an explicit Runge-Kutta discontinuous Galerkin method to achieve sixth-order accuracy when approximating 
the 2D Euler equations on a mesh moved according to a smooth predefined mapping in a fixed
multidimensional domain.

 In subsequent years, high-order approximations in space and time have been proposed for moving
meshes using a range of methods, \emph{e.g.}\ conforming finite elements \cite{ADFR15,ER21,LMQ25,LXY23,RWX25}, finite volumes
\cite{BD14,DBSR17,Gaburro21,GBCKSD20}, discontinuous Galerkin
\cite{BCK20,FHWK21,Fu20,FX22,Gaburro21,GBCKSD20,HR20,HLQZ23,LHQ19,PBP09,PD21,SKBG20,TXL22,ZCHQ20,ZZYSS24,ZX22},
virtual elements \cite{Wells23} and finite differences \cite{LXZ21}.
 The applications in these papers are wide-ranging:
the diffusion equation \cite{LXY23}, advection-diffusion \cite{ER21,LXZ21,Wells23,ZX22},
nonlocal diffusion \cite{ADFR15}, inviscid Burgers' equation \cite{TXL22}, shallow water flows \cite{HLQZ23},
Euler equations \cite{BCK20,BD14,DBSR17,FX22,Gaburro21,GBCKSD20,LHQ19,PBP09,PD21,SKBG20,ZZYSS24},
Stokes equations \cite{RWX25}, incompressible Navier-Stokes equations \cite{FHWK21,Fu20,HR20,LMQ25},
ideal magneto-hydrodynamics \cite{Gaburro21,GBCKSD20} and radiation transfer \cite{ZCHQ20}.
 All, apart from \cite{HR20}, choose an ALE formulation.
 Many restrict themselves to moving meshes within fixed domains, but higher than second-order
accuracy is demonstrated in the presence of a moving boundary in
\cite{ADFR15,ER21,FHWK21,Fu20,HR20,LMQ25,LXY23,LXZ21,LHQ19,PBP09,RWX25,Wells23}. In almost all cases, the boundary
movement is generated by a predefined, smooth, time-dependent mapping between the evolving
physical domain and a reference domain (typically the initial domain). Only \cite{LMQ25,Wells23} demonstrate
higher than second-order accuracy for tests where the boundary/interface movement depends on the
solution.

 Most of the papers cited in the previous paragraph demonstrate higher orders of accuracy
purely through numerical convergence studies, but mathematical analysis to support these
empirical observations of convergence rates is also starting to emerge, \emph{e.g.}\ for conforming finite elements
with one-dimensional nonlocal diffusion \cite{ADFR15} and (on meshes of simplices) multidimensional advection-diffusion
\cite{LXY23}, Stokes flow \cite{RWX25} and two-phase incompressible Navier-Stokes flow with interface tracking \cite{LMQ25},
and for the discontinuous Galerkin method with one-dimensional advection-diffusion \cite{ZX22}.
Elliott and Ranner present a unifying theory for PDEs posed on evolving surfaces and bulk domains in \cite{ER21}
which includes error analysis of finite element approximations of such systems.

\vspace{\baselineskip}
 In this paper, we consider the numerical approximation of linear and nonlinear PDEs
on time-dependent domains, where the movement
of the boundary is determined by the PDE and its boundary conditions. The work builds
on a moving-mesh finite element method originally proposed in \cite{BHJ05a} which was
designed to both accurately approximate the boundary movement and provide an internal
mesh velocity field.
 It generated a Lagrangian mesh velocity field for
problems with moving boundaries by approximating a distributed conservation principle
using finite elements, then recovered the solution on the updated mesh in a manner
which ensured global conservation. Subsequent developments \cite{BHJ11,BHJM09}
used this algorithm to determine the normal velocity
component of the moving boundary, which was then imposed as a boundary condition
while generating the interior mesh motion and the tangential mesh boundary motion using an
alternative approach. Since the mesh movement was no longer Lagrangian, a conservative ALE
step was introduced to evolve the approximate solution in time.

  No mapping to a reference domain is used in the algorithm's implementation, so
the expense of transforming to and from the reference domain is avoided. There is also no
requirement to define what that transformation should be when it isn't
predefined analytically. However, the analytical and computational convenience
of updating the solution on an unchanging reference mesh is lost. This is particularly important
when $p > 1$ and has until now restricted the method to, at most, second-order accuracy for
the prediction of both the boundary movement and the evolution of the solution in the interior.

 This paper provides a first step in extending this mapping-free approach to higher orders
of accuracy for implicit moving boundary problems. The key observation is that, although a
higher-order approximation is required to determine the boundary movement with the appropriate
accuracy, the higher-order information which this generates about the velocity field in the
interior of the domain can be discarded: the internal velocity is taken to be continuous piecewise linear.
This is analogous to using a piecewise affine mapping between the physical domain and a
reference domain, which has recently been observed to be useful in the analysis of
a mapping-based approach with prescribed smooth boundary motion in \cite{LXY23}.
This has enabled us to produce the first arbitrary-order approximations in both space and time to
implicit moving boundary problems within a conservative ALE framework.
These results were obtained without having to adjust the time-stepping schemes
to satisfy a discrete geometric conservation law.

\vspace{\baselineskip}

 The layout of the paper is as follows.
 Section \ref{sec:algorithm} summarises the steps
required to carry out a single time-step of the algorithm. The core steps are the same as
those presented in \cite{BHJ05a,BHJ11}: the main difference is in the computation of the mesh
velocity, which is modified to allow higher than second-order accuracy.
The full, multidimensional, algorithm is presented because the proposed modifications would
apply to any mesh of flat-faced simplices which is allowed to move inside a fixed, flat-faced, boundary.
However, due to our focus on implicit moving boundary problems the test cases we select
are all one-dimensional, to avoid the presence of curved boundaries.
Section \ref{sec:applications} presents the specific discrete forms of the algorithm for three PDEs
with distinct characteristics, selected to test different aspects of the method: (1) the
two-phase Stefan problem (linear diffusion), in which a sharp internal interface moves in
response to local solution gradients; (2) the Crank-Gupta problem (linear reaction-diffusion),
in which boundary movement is governed implicitly by the interaction of an absorption term
with the boundary conditions; (3) the porous medium equation (nonlinear diffusion),
in which the boundary movement is governed by the local solution gradient.
Numerical results are provided which demonstrate, for all three PDEs, up to fifth-order accuracy
in space and time for the approximations of both the solution and the boundary/interface position.
Finally, Section \ref{sec:conc} includes a discussion of points raised in the paper.

\section{The Moving-Mesh Algorithm}
\label{sec:algorithm}

 In this section we summarise the steps of the velocity-based moving-mesh framework
\cite{BHJ05a,BHJ11} on which we build, approximating the general evolutionary PDE,
\begin{equation}
 \frac{\partial u}{\partial t} \;=\; \mathcal{L} u \, ,
 \label{eq:general}
\end{equation}
in which $t \in (t^0,T]$ indicates time and $\mathcal{L}$ a general spatial differential operator
in $\mathbb{R}^d$. The time-dependent support of the evolving solution of Equation \eqref{eq:general}
will be denoted by $\Omega(t) \subset \mathbb{R}^d$, with boundary $\partial \Omega(t)$,
which is allowed to move. A moving coordinate system is defined by
$\mathbf{x} = \mathbf{x}(\mbox{\boldmath$\xi$},t)$, where
$\mathbf{x}(\mbox{\boldmath$\xi$},0) = \mbox{\boldmath$\xi$}$ is the coordinate system
on a fixed reference domain (often taken to be $\Omega(t^0)$).
The time-dependent evolution of $\Omega(t)$ is determined by the velocity field
$\mathbf{v} = \frac{\partial \mathbf{x}}{\partial t}$ and it is this velocity which will be used
in the definition of the moving-mesh finite element method: the reference domain is not required in the
implementation. Initial conditions $u(\mathbf{x},t^0)$ are provided on $\Omega(t^0)$ and
boundary conditions specific to each test will be defined on $\partial \Omega(t)$ for $t \in (t^0,T]$
in Section \ref{sec:applications}.

 The foundations of the algorithm are described in detail in previous publications \cite{BHJ05a,BHJ11,BHJM09},
but we highlight one feature here. Both the mesh velocity and the evolution of the approximate
solution to Equation \eqref{eq:general} are derived from either the finite element formulation of
the conservative ALE equation \cite{BCR75},
\begin{equation} \label{eq:baines_form}
 \frac{\mathrm{d}}{\mathrm{d}t} \int_{\Omega(t)} w \, u \; \mathrm{d}\mathbf{x}
     \;=\; \int_{\Omega(t)} w \left( \mathcal{L} u + \nabla \cdot ( u \mathbf{v} ) \right) \; \mathrm{d}\mathbf{x} \, ,
\end{equation}
or the discrete geometric conservation law \cite{BG04,FGG01,LF96},
\begin{equation} \label{eq:constant_monitor}
 \frac{\mathrm{d}}{\mathrm{d}t} \int_{\Omega(t)} w \; \mathrm{d}\mathbf{x}
     \;=\; \int_{\Omega(t)} w \, \nabla \cdot \mathbf{v} \; \mathrm{d}\mathbf{x} \, ,
\end{equation}
obtained from Equation \eqref{eq:baines_form} when $\mathcal{L}u \equiv 0$ and $u$ is constant in space.

 Equation \eqref{eq:baines_form} is widely used to evolve $u$ on a moving mesh (see, for example,
\cite{ER21,EGP09,FN04}), but a distinguishing feature of the framework we are building on \cite{BHJ05a,BHJ11} is that we
also use Equations \eqref{eq:baines_form} and \eqref{eq:constant_monitor}
to derive the moving-mesh velocity field $\mathbf{v}$. This is done by assuming that the discrete
test functions form a partition of unity and the distribution of the functionals,
\begin{equation} \label{eq:functionals}
 \mu_u(t;w) \;:=\; \int_{\Omega(t)} w \, u \; \mathrm{d}\mathbf{x} \qquad\mathrm{or}\qquad
 \mu_1(t;w) \;:=\; \int_{\Omega(t)} w \; \mathrm{d}\mathbf{x} \, ,
\end{equation}
between the test functions remains constant in time \cite{BHJ05a,BHJM09}. The approximate
velocity field can be generated using polynomials of arbitrarily high degree, but it
cannot be substituted straightforwardly back in to Equation \eqref{eq:baines_form} because
that equation has been derived by assuming
\begin{equation} \label{eq:test_movement}
  \mathcal{D}_t w \;=\; \frac{\partial w}{\partial t} + \mathbf{v} \cdot \nabla w \;=\; 0 \, ,
\end{equation}
\emph{i.e.}\ the test functions move with velocity $\mathbf{v}$
(see, for example, \cite{BG04,FN99,FN04,Gastaldi01}). For general $\mathbf{v}$,
the evolved test functions do not form a basis of a useful space and therefore cannot be
used within a mapping-free finite element method.

 The main contribution of this paper is a method which demonstrates arbitrary-order accuracy in both
space and time. This is enabled by noting that the high-order velocity information is only required
at the domain boundary, so we replace the high-order approximation to $\mathbf{v}$ in the
interior with its piecewise linear interpolant on the current mesh. This corresponds to an
affine mapping between each physical element and its reference element so, for any given
polynomial degree, the basis of the spatial finite element test space on the mesh at time $t$ is mapped to
the basis of the spatial finite element test space of the same polynomial degree on the mesh at time
$t + \Delta t$. Equation \eqref{eq:test_movement} is therefore satisfied by the standard choices of
test functions defined on the evolved mesh at each discrete time level and
integration can be carried out in the physical domain. This modification
can be applied generally in one dimension and on any mesh of simplices in multiple dimensions.
In this paper we focus on results in 1D with a moving boundary: in higher dimensions it is
most likely to be useful when the mesh is moved inside a domain with a fixed boundary because allowing free movement would generate curved boundaries even where they were initially flat.

\vspace{\baselineskip}

 We now summarise the algorithm for a single time-step, using a subscript $\cdot_h$ to indicate
spatial approximations on a computational domain $\Omega_h(t) \approx \Omega(t)$, consisting of a moving tessellation
of simplices (line segments in one dimension) with boundary $\partial \Omega_h(t) \approx \partial \Omega(t)$.
In the description of the algorithm below we will suppress the dependence on time because within each step the spatial
integrals are all evaluated at the specific, discrete, time levels generated by the application of a finite difference approximation
in time. The description below outlines the multidimensional framework introduced in \cite{BHJ05a,BHJ11,BHJM09}.
This framework has also been used to produce a second-order algorithm on polygonal meshes with the
virtual element method \cite{WHC24}, but is enhanced in this work by introducing $\tilde{\mathbf{v}}$,
the piecewise linear interpolant on the computational of the mesh velocity field $\mathbf{v}$, in the evaluation of $u$
at the new time level, the extra step required for high-order accuracy.

 The evolution of the mesh and solution can be expressed, using a method of lines approach, in the form of a system
of ordinary differential equations in time \cite{BHJ05a,BHJ11,WHC24}:
\begin{equation}
 \frac{\mathrm{d}}{\mathrm{d}t} \left( \begin{array}{c} \mathbf{x}_h \\ \bm{\mu}_{h,u} \end{array} \right)
          \;=\; \mathbf{F}(\mathbf{x}_h,\bm{\mu}_{h,u}) \, ,
 \label{eq:ODEupdate}
\end{equation}
in which $\mathbf{x}_h$ contains the degrees of freedom which define the moving coordinate system (the mesh node
positions in our case) and $\bm{\mu}_{h,u}$ is a vector of discrete values (of the first integral defined by Equation \eqref{eq:functionals})
defined by the spatial basis functions $w_h(\mathbf{x})$ of the finite-dimensional test space at the current time.
The right-hand side function $\mathbf{F}(\mathbf{x}_h,\bm{\mu}_{h,u})$ is computed via the following steps.

\begin{enumerate}

 \item Given $\mathbf{x}_h$ (a coordinate system defined by the current positions of the mesh nodes)
  and $\bm{\mu}_{h,u}$ (the distribution of the `mass', $\theta_{h,u} = \int_{\Omega_h} u_h \, \mathrm{d}\mathbf{x}$,
  between the discrete test functions), recover $u_h$ (the approximate solution on the current mesh) by solving
  \begin{equation}
    \int_{\Omega_h} w_h \, u_h \; \mathrm{d}\mathbf{x} \;=\; \mu_{h,u}(w_h) \, .
    \label{eq:recoversolution_discrete}
  \end{equation}
  This step requires that we keep track of the evolution of the local masses (in Step 5), even when they are not being used in Step 3 to determine the mesh movement.
  

 \item Given $u_h$, recompute the distribution of the integral of a monitor function $\mathbb{M}(u_h)$ between
  the test functions on the current mesh, using
  \begin{equation}
   c_{h,\mathbb{M}}(w_h) \;=\; \int_{\Omega_h} w_h \, \mathbb{M}(u_h) \; \mathrm{d}\mathbf{x}
              \left/ \int_{\Omega_h} \mathbb{M}(u_h) \; \mathrm{d}\mathbf{x} \right.
                                       \;:=\; \frac{\mu_{h,\mathbb{M}}(w_h)}{\theta_{h,\mathbb{M}}} \, .
  \end{equation}
  In this work we consider two cases with different advantages.

  $\mathbb{M}(u) = u$, which we refer to as the mass monitor, approximates the Lagrangian `flow' velocity of the
  PDE when $\mathcal{L} = \nabla \cdot \mathbf{f}$ for some flux $\mathbf{f}$ in Equation \eqref{eq:general}
  \cite{BHJ05a,BHJ11,BHJM09}.
  This choice has two benefits: (i) it allows us to predict the movement
  of free boundaries where there is no kinematic condition explicitly relating the boundary velocity to the local
  solution; (ii) it reduces interpolation error between time-steps in the discretisation because the
  mesh (and hence the solution) is transported with the velocity field inherent to the PDE.

  $\mathbb{M}(u) = 1$, which we refer to as the area monitor, approximately preserves the distribution of mesh
  element length/area/volume when boundary velocities are provided \cite{MHJ11}.
  It cannot be used to predict boundary velocities but can be used to generate internal velocities when an explicit
  form is provided for the boundary velocity or where the boundary movement has already been computed using
  $\mathbb{M}(u) = u$.

  This step is omitted when $\dot{\theta}_{h,\mathbb{M}} = 0$ because there is no longer a
  requirement for $c_{h,\mathbb{M}}(w_h)$ in the next steps.

 \item This step is only applied when no kinematic boundary condition (\emph{i.e.}\ one which
  explicitly includes the boundary velocity $\mathbf{v}$) is available to provide $u_h \, \mathbf{v}_h \cdot \mathbf{n}$
  on the moving boundary, where $\mathbf{n}$ is the spatial outward unit normal to $\partial \Omega_h$.
  It computes a piecewise linear interpolant $\tilde{\mathbf{v}}_h$ of a higher-degree velocity field
  to approximate the physical boundary
  motion and provide mesh velocities for the interior (although the latter may be replaced in Step 4). It uses
  Equation \eqref{eq:baines_form} and assumes that the velocity field is curl-free, so it
  can be written as $\mathbf{v} = \nabla \phi$ for some velocity potential $\phi$ \cite{BHJ05a}.

  (a) Given $u_h$, determine $\phi_h$ (the discrete mesh velocity potential) by solving
  \begin{eqnarray}
    c_{h,u}(w_h) \, \dot{\theta}_{h,u} \,+\, \int_{\Omega_h} u_h \, \nabla w_h \cdot \nabla \phi_h \; \mathrm{d}\mathbf{x}
       &=& \int_{\Omega_h} w_h \, ( \mathcal{L} u )_h \; \mathrm{d}\mathbf{x}
      \,+\, \int_{\partial \Omega_h} w_h \, u_h \, \mathbf{v}_h \cdot \mathbf{n} \; \mathrm{d}s \, , \nonumber \\
    \dot{\theta}_{h,u} &=& \int_{\Omega_h} ( \mathcal{L} u )_h \; \mathrm{d}\mathbf{x} \,+\,
        \int_{\partial \Omega_h} u_h \, \mathbf{v}_h \cdot \mathbf{n} \; \mathrm{d}s \, , \label{eq:potentialequation_discrete}
  \end{eqnarray}
  in which $(\mathcal{L} u)_h$ simply indicates that the spatial operator $\mathcal{L} u$ will be approximated: its precise form will be given for each test case in Section \ref{sec:applications}.
  Note that, even when there is no explicit kinematic condition on the moving boundary, the boundary
  integrals above disappear due to the conditions on $u_h$ in all cases in Section \ref{sec:applications}.
  The boundary integrals are zero on fixed boundaries.
  

  (b) Given $\phi_h$, recover $\tilde{\mathbf{v}}_h$ (a piecewise linear mesh velocity) by first solving
  \begin{equation}
   \int_{\Omega_h} z_h \, \mathbf{v}_h \; \mathrm{d}\mathbf{x}
       \;=\; \int_{\Omega_h} \, z_h \, \nabla \phi_h \; \mathrm{d}\mathbf{x} \, ,
   \label{eq:velocityrecovery_discrete}
  \end{equation}
  and then taking $\tilde{\mathbf{v}}_h$ to be the continuous piecewise linear interpolant of $\mathbf{v}_h$
  on the current mesh. This gives boundary velocities (which can be imposed strongly when known) 
  used as the right-hand side for $\dot{\mathbf{x}}_h$ in Equation \eqref{eq:ODEupdate}. This step
  is valid in more than one dimension if the spatial mesh is composed of simplices.
  Note that if, instead of interpolating, $\tilde{\mathbf{v}}_h$ is recovered by
  projection of $\nabla \phi_h$ directly on to the space of piecewise linear polynomials on the
  current mesh, the results are less accurate and demonstrate lower orders of convergence.
  


 \item This step is applied when $\mathbf{v}_h \cdot \mathbf{n}$ is provided on $\partial \Omega_h$ by the boundary conditions or
  $\tilde{\mathbf{v}}_h \cdot \mathbf{n}$ has been computed on the boundary in Step 3, but a different procedure
  is preferred for computing the interior velocity field.

  (a) Given $u_h$ and $\tilde{\mathbf{v}}_h \cdot \mathbf{n}$ on $\partial \Omega_h$, determine
  $\phi_h$ (a discrete mesh velocity potential for interior motion)
  by either (i) solving Equations \eqref{eq:potentialequation_discrete} if $\mathbb{M}(u) = u$ is
  used to drive the internal mesh movement or, (ii) if $\mathbb{M}(u) = 1$ is used to drive the internal mesh movement,
  solving
  \begin{eqnarray}
    c_{h,1}(w_h) \, \dot{\theta}_{h,1} \,+\, \int_{\Omega_h} \nabla w_h \cdot \nabla \phi_h \; \mathrm{d}\mathbf{x}
      &=& \int_{\partial \Omega_h} w_h \, \tilde{\mathbf{v}}_h \cdot \mathbf{n} \; \mathrm{d}s \, , \nonumber \\
    \dot{\theta}_{h,1} &=& \int_{\partial \Omega_h} \tilde{\mathbf{v}}_h \cdot \mathbf{n} \; \mathrm{d}s \, .
   \label{eq:potentialequation2_discrete}
  \end{eqnarray}

  (b) Given $\phi_h$, recover $\tilde{\mathbf{v}}_h$ (a piecewise linear mesh velocity for interior motion)
  as in step 3(b). This gives internal velocities which
  can be used as the right-hand side for $\dot{\mathbf{x}}_h$ in Equation \eqref{eq:ODEupdate}.

 \item[(5)] Given $u_h$ and $\tilde{\mathbf{v}}_h$, compute $\dot{\mu}_{h,u}(w_h)$ (the rate of change
  of `mass' associated with each discrete test function) using
  \begin{equation}
    \dot{\mu}_{h,u}(w_h) \;=\; \int_{\Omega_h} w_h \, ( \mathcal{L} u )_h \; \mathrm{d}\mathbf{x}
       \,-\, \int_{\Omega_h} u_h \, \nabla w_h \cdot \tilde{\mathbf{v}}_h \; \mathrm{d}\mathbf{x}
       \,+\, \int_{\partial \Omega_h} w_h \, u_h \, \tilde{\mathbf{v}}_h \cdot \mathbf{n} \; \mathrm{d}s \, .
    \label{eq:ALEupdate_discrete}
  \end{equation}

\end{enumerate}

 The system of ODEs given by Equation \eqref{eq:ODEupdate} can be integrated in time using the most
appropriate time-stepping scheme for the problem. In the following work Strong-Stability-Preserving Runge-Kutta
(SSP-RK) schemes \cite{GST01} are chosen, where the order of accuracy is selected
so that the error in the temporal discretisation reduces at least as fast as the error in the spatial
discretisation for refinement paths which are chosen to ensure stability as the mesh size tends to zero.
These time-stepping schemes are not guaranteed to satisfy a discrete geometric conservation law.
They could be adjusted in the manner of \cite{BG04} to do so, but the numerical experiments in
Section \ref{sec:applications} demonstrate that this is not required for retaining the expected orders
of accuracy in space and time. Standard Lagrange basis functions are chosen for the spatial discretisation
because they define a partition of unity and link each test function with a distinct spatial position.

 We further note that any mesh movement strategy could be used to redistribute the internal mesh
nodes in step 4 above, {\emph e.g.}\ to adapt to or equidistribute a local error indicator or to
regularise the mesh movement. The new node positions would then be used to generate nodal velocities
by assuming that these velocities were constant within each time-step, and $\tilde{\mathbf{v}}_h$
would be their continuous piecewise linear interpolant. However, the aim of this paper is to demonstrate high-order
accuracy and the optimisation of the node positions is a future task.

\section{Applications}
\label{sec:applications}

 We consider an approximation $\Omega_h(t)$ to $\Omega(t)$ with boundary
$\partial \Omega_h(t) \approx \partial \Omega(t)$, divided into a fixed part,
$\partial \Omega_{h,F} \approx \partial \Omega_F$ (which are equal in one dimension), and a moving part,
$\partial \Omega_{h,M}(t) \approx \partial \Omega_M(t)$, which satisfy
\begin{equation}
 \partial \Omega_h(t) \;=\; \partial \Omega_{h,M}(t) \cup \partial \Omega_{h,F} \qquad \mbox{and} \qquad
 \partial \Omega_{h,M}(t) \cap \partial \Omega_{h,F} \;=\; \emptyset \, .
\end{equation}
The test cases considered in this paper are all one-dimensional so the shape of the boundary is not approximated.
However, the notation above is retained to indicate that the boundary position is approximate.
The details of each discrete form are presented below for completeness and to highlight where the
linear interpolant of the mesh velocity field is used.

 We present space-time errors for both the solution and the boundary position.
To account for the use of finite differences in time and finite elements in space, the errors
are measured using mixed discrete/continuous $L^2$-norms of the form
\begin{eqnarray}
  \mbox{Error}_u &=& \left( \sum_{n=1}^{N_t} \Delta t \int_{\Omega_h(t^n)} [ u(\mathbf{x},t^n) - u_h(\mathbf{x},t^n) ]^2 \; \mathrm{d}\mathbf{x} \right)^{\frac{1}{2}} \, , \nonumber \\
  \mbox{Error}_x &=& \left( \sum_{n=1}^{N_t} \Delta t \, | \mathbf{x}_{\partial}(t^n) - \mathbf{x}_{\partial,h}(t^n) |^2 \right)^{\frac{1}{2}} \, , 
  \label{eq:analyticerror}
\end{eqnarray}
in which $\Delta t$ is the time-step length (chosen to be a fixed constant during each
simulation when estimating convergence rates), $N_t$ is the number of
time-steps taken, $u$ and $u_h$ are the exact and approximate solutions to the PDE and
$\mathbf{x}_{\partial}$ and $\mathbf{x}_{\partial,h}$ are the exact and approximate boundary positions
(single points in one dimension).
We have also monitored the convergence behaviour of the $L^2$-norm errors sampled at the final
times of the simulations (not presented here) and note that it matches that of the space-time errors.

 In all simulations we apply an explicit time-stepping method to a
second-order parabolic PDE problem so, as the spatial mesh size $h$ is reduced, the
time-step is reduced according to $\Delta t \propto h^2$, each time the element lengths in the
initial mesh are refined. This relationship between $\Delta t$ and $h$ means that it is possible
to generate the mesh convergence results using time-stepping schemes
of lower order than the spatial discretisation. Forward Euler (which is first-order) was enough
for the piecewise linear spatial approximation ($p = 1$), Heun's method (second-order) was
used for $p = 2,3$ and a third-order strong stability-preserving Runge-Kutta method \cite{GST01}
was used for $p = 4$.
 For a finite element approximation using polynomial basis functions of degree $p$
the approximation was initialised with the continuous piecewise interpolant of the analytical
solution at $p+1$ uniformly distributed points spanning each element.
 Unless stated otherwise, in the mesh convergence studies the spatial meshes are refined
uniformly, with the initial value of $h$ repeatedly divided by 2.

\vspace{\baselineskip}
  Three applications are chosen to illustrate different aspects of the algorithm.
All have suitable analytic solutions with compact, time-dependent support which can be used
to test the method.

\subsection{Change of Phase}

  A two-phase Stefan problem, which models change of phase between liquid and solid, is a linear PDE
with a kinematic boundary condition at the moving interface relating the interface velocity to
the gradient of the solution on either side. The interface condition provides an explicit expression
for the normal velocity component of the moving interface, which is used directly in Equations
\eqref{eq:potentialequation2_discrete}. Unlike in \cite{BHJM09}, we use this test to investigate
the method's behaviour when $\mathbb{M}(u) = 1$ is used to drive
the internal mesh movement, so each mesh element retains, approximately, its length relative
to all of the other elements (and the domain size). The initial mesh is locally refined close to the
moving boundary and this choice of $\mathbb{M}(u)$ ensures that the mesh retains this property as it evolves.

\subsubsection{The Stefan Problem}

 Consider a domain divided into regions of solid (S) or liquid (L) phase, denoted respectively by
$\Omega_S(t)$ and $\Omega_L(t)$. The two-phase Stefan problem (described in detail in
\cite{Cra84}) models transition between liquid and solid phases
across an interface, denoted here by $\partial \Omega_M(t)$, in
the interior of the problem domain which moves as time progresses. It is modelled by \cite{BCFP73}
\begin{eqnarray}
 K_S \, u_t &=& \nabla \cdot ( k_S \, \nabla u )
          \hspace{10mm} {\mbox{in}} \;\;\; \Omega_S(t) \nonumber \\
 K_L \, u_t &=& \nabla \cdot ( k_L \, \nabla u )
          \hspace{10mm} {\mbox{in}} \;\;\; \Omega_L(t) \, ,
 \label{eq:Stefan}
\end{eqnarray}
in which $K_S,K_L$ represent the volumetric heat capacities, $k_S,k_L$ the
thermal conductivities and $u$ the temperature. The interface conditions are
$u = u_M$ and
\begin{equation}
 k_S \, \nabla u_S \cdot \mathbf{n} - k_L \, \nabla u_L \cdot \mathbf{n}
   \;=\; \lambda \, \mathbf{v} \cdot \mathbf{n} \, ,
 \label{eq:StefanBC}
\end{equation}
where $\mathbf{n} = \mathbf{n}(t)$ is a unit normal to the moving interface, $\lambda$ is
the heat of phase-change per unit volume, and $\mathbf{v} \cdot \mathbf{n}$ is the
normal velocity of the interface. All parameters are assumed here to be
constant within their respective phases
and subscripts are used in Equation \eqref{eq:StefanBC} to indicate the phase in which the gradient is evaluated.
The fixed (external) boundaries for each phase are denoted by $\partial \Omega_{F_S}$ and $\partial \Omega_{F_L}$.

\subsubsection{Algorithm for the Stefan Problem}

 In this work we propose an alternative to \cite{BHJM09}, in which the algorithm was equivalent to using $\mathbb{M}(u) = u + K$, for some large enough positive constant $K$, to determine the mesh movement. Instead, we select
$\mathbb{M}(u) = 1$, a simpler and more robust choice which preserves element
lengths and, unlike in \cite{BHJM09}, allows us to impose $u = u_M = 0$ at the moving interface without
adjusting the algorithm. Hence, when the initial mesh is adapted locally to resolve the interface between the phases, it
remains adapted in a way which resolves the interface as it moves.
With $\mathbb{M}(u) = 1$, the velocity potential is found using Equations \eqref{eq:potentialequation2_discrete}
which, when combined with the Stefan condition \eqref{eq:StefanBC}, give new equations for the velocity potential:
\begin{eqnarray}
 & & \hspace*{-30mm} c_{h,S,1}(w_h) \, \dot{\theta}_{h,S,1}(t)
  \;+\;  \int_{\Omega_{h,S}(t)} \nabla \phi_h \cdot \nabla w_h \; \mathrm{d}\mathbf{x} \\
 &=& \int_{\partial \Omega_{h,M}(t)} \frac{w_h}{\lambda}
    \left\{ k_S ( \nabla u_h )_S \cdot \mathbf{n}_S - k_L ( \nabla u_h )_L \cdot \mathbf{n}_S \right\} \; \mathrm{d}s \, , \nonumber \\
 & & \hspace*{-30mm} c_{h,L,1}(w_h) \, \dot{\theta}_{h,L,1}(t)
  \;+\;  \int_{\Omega_{h,L}(t)} \nabla \phi_h \cdot \nabla w_h \; \mathrm{d}\mathbf{x} \\
 &=& \int_{\partial \Omega_{h,M}(t)} \frac{w_h}{\lambda}
    \left\{ k_S ( \nabla u_h )_S \cdot \mathbf{n}_L - k_L ( \nabla u_h )_L \cdot \mathbf{n}_L \right\} \; \mathrm{d}s \, , \nonumber
 \label{eq:stefan_potential}
\end{eqnarray}
along with
\begin{eqnarray}
 \dot{\theta}_{h,S,1}(t) &=& \int_{\partial \Omega_{h,M}(t)} \frac{1}{\lambda}
    \left\{ k_S ( \nabla u_h )_S \cdot \mathbf{n}_S - k_L ( \nabla u_h )_L \cdot \mathbf{n}_S \right\} \; \mathrm{d}s \, , \nonumber \\
 \dot{\theta}_{h,L,1}(t) &=& \int_{\partial \Omega_{h,M}(t)} \frac{1}{\lambda}
    \left\{ k_S ( \nabla u_h )_S \cdot \mathbf{n}_L - k_L ( \nabla u_h )_L \cdot \mathbf{n}_L \right\} \; \mathrm{d}s \, .
 \label{eq:stefan_theta1dot}
\end{eqnarray}
The vectors $\mathbf{n}_S$ and $\mathbf{n}_L$ are the unit normal vectors at the
moving interface, pointing outwards from their respective domains. Computing the
velocity potential $\phi_h$ separately in the different regions is essential because
it allows the meshes in the solid and liquid phases to independently retain their
initial distributions of element lengths even though the relative sizes of the
solid and liquid domains change with time.

 As described in \cite{BHJM09}, the mesh velocities are recovered by solving a single system. The new
element here is the replacement of the solution to \eqref{eq:velocityrecovery_discrete} in the interiors of the liquid and
solid domains by its continuous piecewise linear interpolant. The system is constrained by
$\tilde{\mathbf{v}}_h = \mathbf{0}$ at each (fixed) end of the interval and, at the moving interface,
\begin{equation}
    \lambda \int_{\partial \Omega_{h,M}(t)} w_h \, \tilde{\mathbf{v}}_h \cdot \mathbf{n} \; \mathrm{d}s
  \;=\; k_S \int_{\partial \Omega_{h,M}(t)} w_h \, (\nabla u_h)_S \cdot \mathbf{n} \; \mathrm{d}s
  \;-\; k_L \int_{\partial \Omega_{h,M}(t)} w_h \, (\nabla u_h)_L \cdot \mathbf{n} \; \mathrm{d}s \, ,
 \label{eq:xdotbcsvec}
\end{equation}
where $\mathbf{n}$ is a normal vector to the moving interface.

 The dependent variable $u_h$ is recovered separately in each region of the domain,
$\Omega_{h,S}(t)$ and $\Omega_{h,L}(t)$, since they are decoupled by the Dirichlet boundary
condition, $u_h = u_M$, imposed on the interface \cite{BHJM09}.
If $\mu_{h,S,u}$ and $\mu_{h,L,u}$ are the local masses in the solid and liquid phases, respectively,
then applying the boundary conditions $u_h = u_M = 0$ on $\partial \Omega_M(t)$ and
$\tilde{\mathbf{v}}_h \cdot \mathbf{n} = 0$ on $\partial \Omega_{F_S} \cup \partial \Omega_{F_L}$ leads,
via Equation \eqref{eq:ALEupdate_discrete} and an interpolated mesh velocity, to
\begin{eqnarray}
 \dot{\mu}_{h,S,u}(t;\hat{w}_h)
    &=& \int_{\partial \Omega_{h,M}(t)} \hat{w}_h \, \kappa_S (\nabla u_h)_S \cdot \mathbf{n}_S \; \mathrm{d}s
         \,+\, \int_{\partial \Omega_{h,F}} \hat{w}_h \, \kappa_S  \nabla u_h \cdot \mathbf{n} \; \mathrm{d}s \nonumber \\
    & &  \hspace*{15mm} \,-\, \int_{\Omega_{h,S}(t)} \nabla \hat{w}_h \cdot ( \kappa_S \nabla u_h + u_h \, \tilde{\mathbf{v}}_h ) \; \mathrm{d}\mathbf{x} \, , \\
 \dot{\mu}_{h,L,u}(t;\hat{w}_h)
    &=& \int_{\partial \Omega_{h,M}(t)} \hat{w}_h \, \kappa_L (\nabla u_h)_L \cdot \mathbf{n}_L \; \mathrm{d}s
         \,+\, \int_{\partial \Omega_{h,F}} \hat{w}_h \, \kappa_L  \nabla u_h \cdot \mathbf{n} \; \mathrm{d}s \nonumber \\
    & &  \hspace*{15mm} \,-\, \int_{\Omega_{h,L}(t)} \nabla \hat{w}_h \cdot ( \kappa_L \nabla u_h + u_h \, \tilde{\mathbf{v}}_h ) \; \mathrm{d}\mathbf{x} \, ,
 \label{eq:Stefan:ALE}
\end{eqnarray}
where $\kappa_* = k_* / K_*$ for $* \in \{S,L\}$ and $(\nabla u_h)_S$, $(\nabla u_h)_L$ are the gradients of
the approximate solution on, respectively, the solid and liquid sides of the moving interface.
In the derivation of these equations $(\mathcal{L} u)_h$ has been defined by the discrete form
\begin{equation}
 \int_{\Omega_{h,*}} \hat{w}_h \, (\mathcal{L} u)_h \; \mathrm{d}\mathbf{x}
    \;=\; \int_{\partial \Omega_{h,*}} \hat{w}_h \, \kappa_* ( \nabla u_h )_* \cdot \mathbf{n}_* \; \mathrm{d}s
    \;-\; \int_{\Omega_{h,*}} \nabla \hat{w}_h \cdot ( \kappa_* \nabla u_h )\; \mathrm{d}\mathbf{x} \, .
\end{equation}

We choose $\hat{w}_h$ to be the modified test functions of \cite{HBJ09}, so that
$u$ is specified strongly on Dirichlet boundaries without losing conservation. The same
convergence rates are seen in the numerical experiments when the mesh is refined if standard
test functions are used with strongly imposed Dirichlet boundary conditions, at the expense
of global conservation. The positions of the mesh nodes
and the values of $\mu_{h,S,u}(t;\hat{w}_h)$ and $\mu_{h,L,u}(t;\hat{w}_h)$ are then evolved to a
later time $t + \Delta t$, after which the dependent variable is recovered from
\begin{eqnarray}
 \int_{\Omega_{h,S}(t+\Delta t)} \hat{w}_h \, u_h \; \mathrm{d}\mathbf{x} &=& \mu_{h,S,u}(t+\Delta t;\hat{w}_h) \, , \nonumber \\
 \int_{\Omega_{h,L}(t+\Delta t)} \hat{w}_h \, u_h \; \mathrm{d}\mathbf{x} &=& \mu_{h,L,u}(t+\Delta t;\hat{w}_h) \, .
 \label{eq:psi_2phase}
\end{eqnarray}

\subsubsection{Numerical Results}

 An exact solution to the one-dimensional two-phase Stefan problem is given
in \cite{CJ59} as
\begin{eqnarray}
 u_S &=& u^* \left( 1 - \frac{er\!f \, ( x / (2 \sqrt{\kappa_S \, t}) )}
                             {er\!f \, \phi} \right) \nonumber \\
 u_L &=& u_0 \left( 1 - \frac{er\!f\!c \, ( x / (2 \sqrt{\kappa_L \, t}) )}
                             {er\!f\!c \, ( \phi \sqrt{\kappa_S / \kappa_L})} \right) \, ,
 \label{eq:exactStefan}
\end{eqnarray}
where $\phi$ is the root of
\begin{equation}
 \frac{e^{-\phi^2}}{er\!f \, \phi}
   + \frac{k_L}{k_S} \sqrt{\frac{\kappa_S}{\kappa_L}}
       \frac{u_0 e^{-\kappa_S \phi^2 / \kappa_L}}{u^* \, er\!f\!c \, ( \phi \sqrt{\kappa_S / \kappa_L} )}
   + \frac{\phi \, \lambda \, \sqrt{\pi}}{K_S \, u^*} \;=\; 0 \, ,
\end{equation}
and $er\!f$ and $er\!f\!c$ are the standard error and complementary error
functions, respectively. The position of the moving interface is given by
\begin{equation}
 s(t) \;=\; 2 \, \phi \, \sqrt{\kappa_S \, t} \, .
\end{equation}
Following previous work \cite{BHJM09,Fur80,MR00}, we select
\begin{equation}
 k_S \;=\; 2.22  \, , \quad k_L \;=\; 0.556 \, , \quad
 K_S \;=\; 1.762 \, , \quad K_L \;=\; 4.226 \, , \quad
 \lambda \;=\; 338 \, ,
\end{equation}
with $u^* = -20$ and $u_0 = 10$, so $\phi \approx 0.205426929376498$.
To avoid the singularity at $t = 0$, we take $t^0 = 0.0012$, at which time the interface is at $x \approx 0.01597539$. 
In the following numerical experiments, the computational
domain is $x \in [0,1]$, a Dirichlet condition ($u = u^* = -20$) is imposed
at the fixed boundary at $x = 0$ and a Neumann condition (the exact value of
$u_x(1,t)$, derived from $u_L$ in Equation \eqref{eq:exactStefan})
is imposed at the fixed boundary at $x = 1$. The Stefan condition \eqref{eq:StefanBC}
is applied with $u = u_M = 0$ at the moving interface.

\begin{figure}[htbp]
 \centerline{\includegraphics[width=0.65\textwidth]{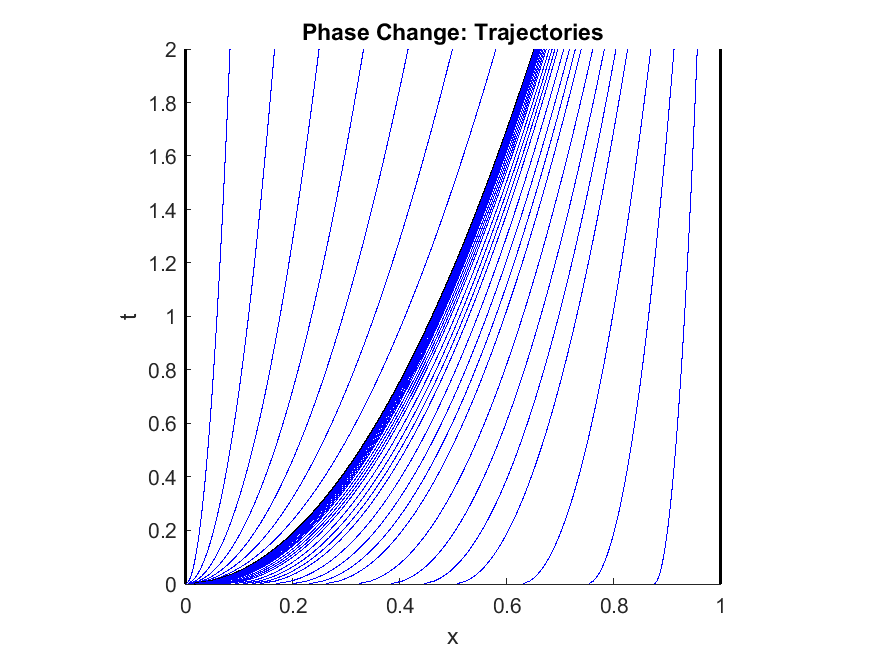} \hspace{-19mm}
          \includegraphics[width=0.65\textwidth]{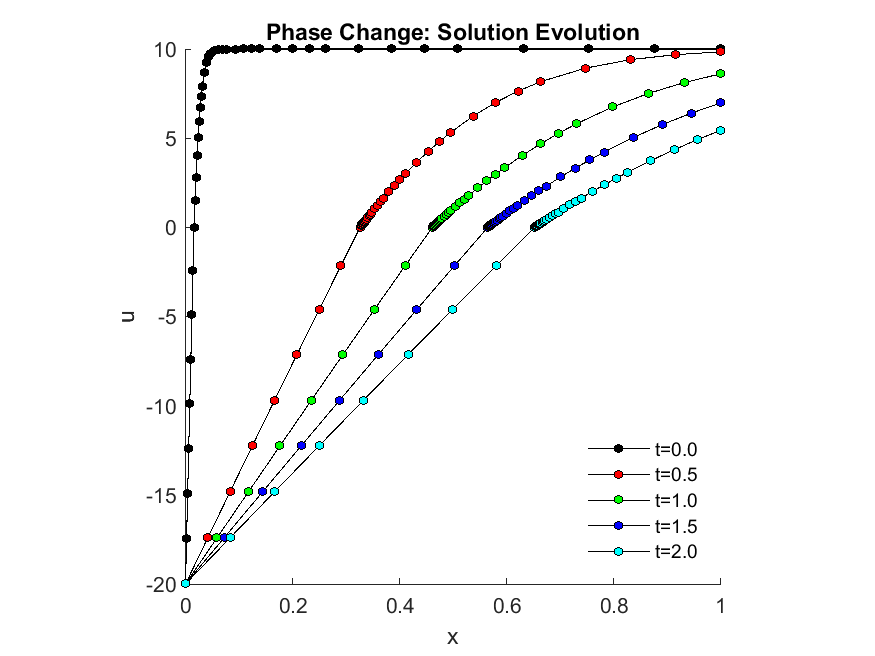}}
 \caption{One-dimensional Stefan problem approximated on a 40-element mesh, initially generated using
          $h$-refinement to the right of the moving interface: node trajectories (left);
          snapshots of solution (right).}
 \label{fig:Stefan_evolution}
\end{figure}

\vspace{\baselineskip}
 The initial meshes are uniform in the solid phase (on the left of the
moving interface) and adaptively refined in the liquid phase (on the right). The coarsest
mesh used has 10 elements, 2 on the left of the interface and 8 on the right.
The 8-element mesh on the right is generated by recursively bisecting the interval adjacent
to the moving interface 7 times, to adapt it locally to the moving interface. All other
initial meshes are uniform refinements of this. The evolution of the node positions and
the solution are shown in Figure \ref{fig:Stefan_evolution} for the 40-element mesh.

\begin{figure}[htbp]
 \centerline{\includegraphics[width=0.65\textwidth]{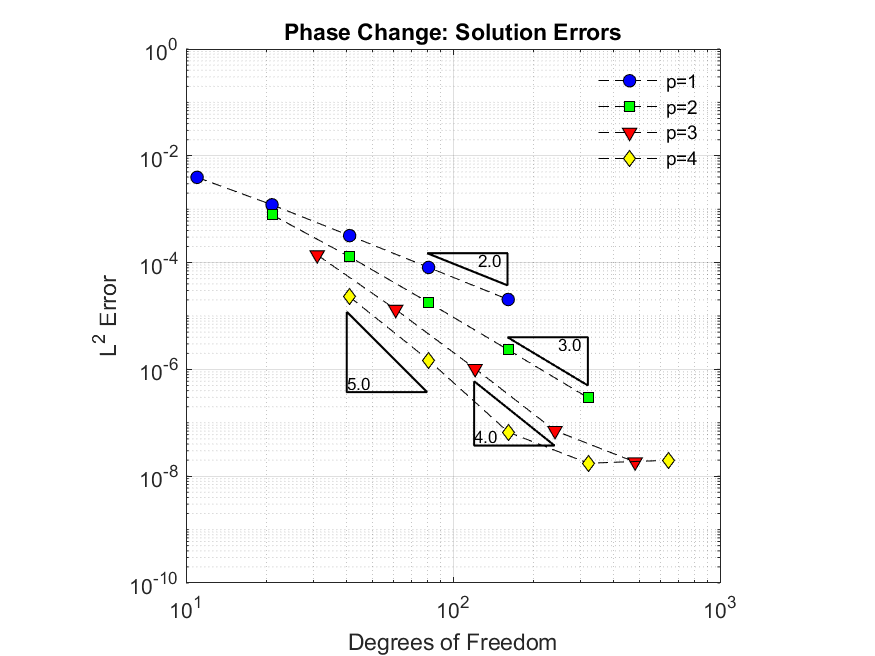} \hspace{-19mm}
          \includegraphics[width=0.65\textwidth]{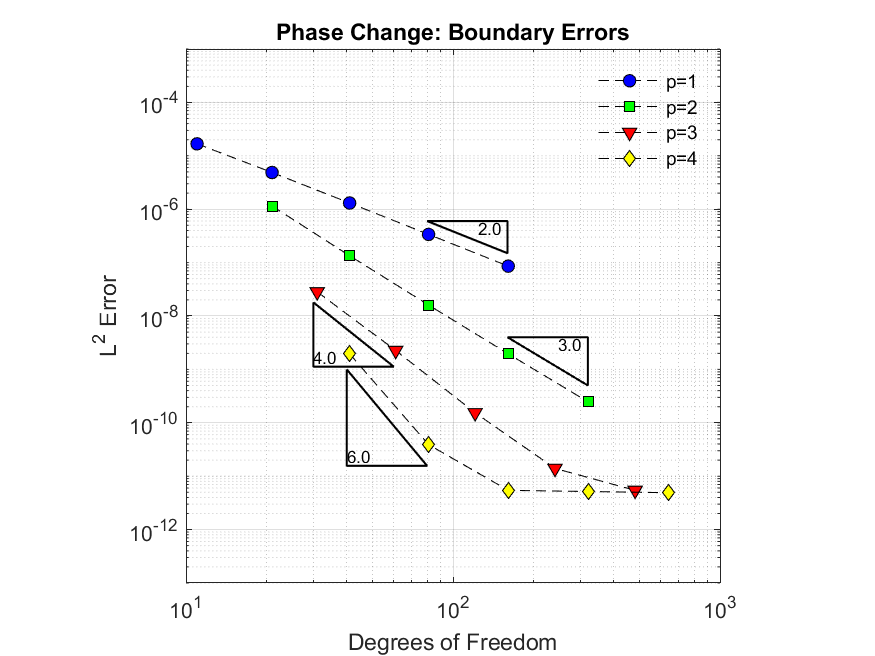}}
 \vspace{5mm}
 \centerline{\includegraphics[width=0.65\textwidth]{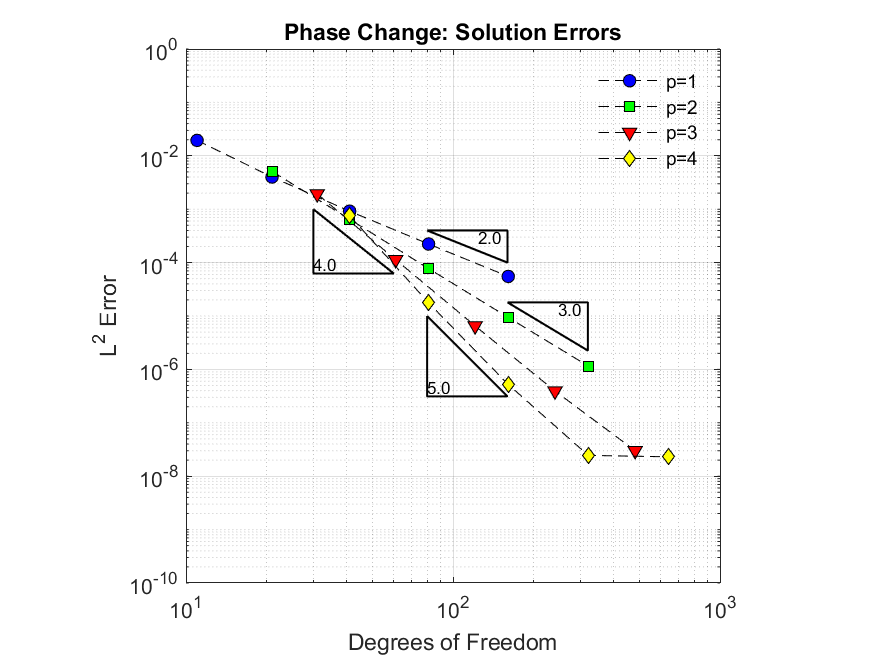} \hspace{-19mm}
          \includegraphics[width=0.65\textwidth]{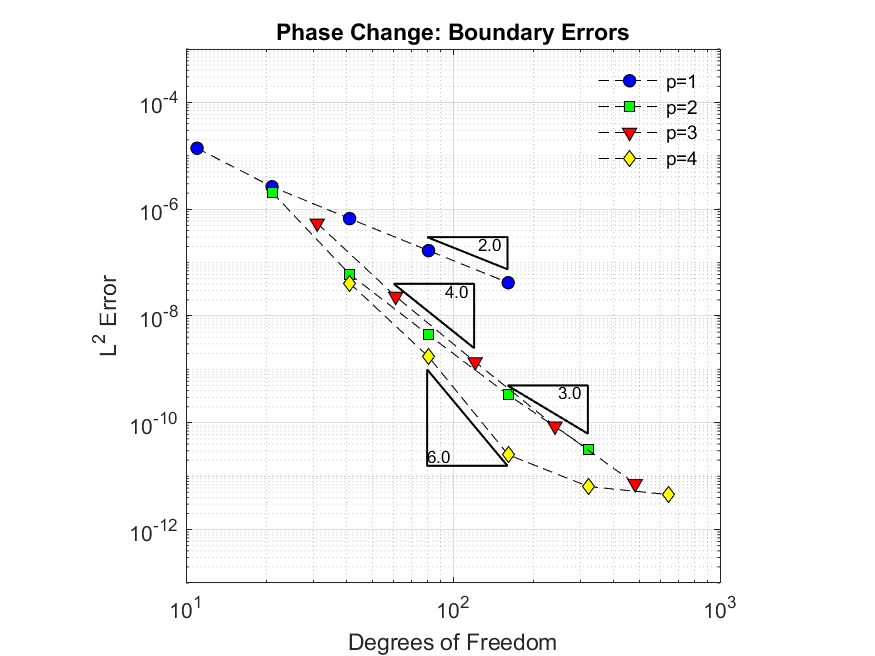}}
 \caption{Mesh convergence plots for the one-dimensional Stefan problem approximated on
          the interval $t \in [t^0,t^0 + 0.01]$: solution error (left); interface position error (right).
          The two graphs at the top show results for initial meshes generated using $h$-refinement,
          and the two at the bottom show results for initial meshes generated using a
          geometric progression.}
 \label{fig:Stefan_errors}
\end{figure}

 This mesh differs from \cite{BHJM09}, in which a uniform mesh was used
on the left of the moving interface and the mesh on the right had the same number of
elements. The two elements on either side of the moving interface were given the same
length, then the remaining element lengths to the right of the initial interface position were generated
as a geometric progression with the common ratio computed to exactly fill the interval between
the interface and the right-hand boundary. Convergence results obtained on these meshes
using the modified algorithm of this paper are also shown in Figure
\ref{fig:Stefan_errors} for comparison.
 Due to very small element size close to the interface, $\Delta t = 8 \times 10^{-7}$
is used for the coarsest (10-element) mesh.

 The graphs in Figure \ref{fig:Stefan_errors} show the expected orders of convergence
for $p = 1,2,3,4$ for both sets of meshes when considering the errors in the solution. 
The same patterns are seen for the errors in the interface position except for $p=4$
which appears to converge at a slightly higher rate. This may be due to
rounding error, which has a prominent effect here because of the very small size
of the time-steps (and hence the very small size of the increments to the solution and
interface position each time a time-step is taken). This is why the error profiles
flatten out at approximately $10^{-8}$ for the solution and $10^{-11}$ for the
interface position. It is also possible that the temporal error is dominating the spatial
error when $p = 4$, since the combination of a third-order time-stepping method and
$\Delta t \propto h^2$ would mimic the convergence behaviour of a sixth-order method.

\subsection{Absorption-Diffusion}

 The Crank-Gupta problem models coupled absorption and diffusion. It provides
a linear PDE with a sink term (so mass is not conserved) and there is no explicit kinematic boundary
condition, so the boundary movement is generated implicitly through the interaction between the
PDE and its boundary conditions. It is used here to demonstrate that the choice of $\mathbb{M}(u) = u$
can be used to accurately determine the normal boundary velocity when it is not provided by a kinematic
condition of the form $\mathbf{v} \cdot \mathbf{n} = \mathcal{B}u$, where $\mathcal{B}$ is some
operator, and in cases where $\dot{\theta}_u \neq 0$.

\subsubsection{The Crank-Gupta Problem}

 The Crank-Gupta problem \cite{BCR75,CG72,Ock96} provides a
system where the total mass changes due to the presence of a source/sink term and is given by
\begin{equation}
 \frac{\partial u}{\partial t} \;=\; \nabla^2 u - 1
   \hspace*{10mm} {\mbox{in}} \;\;\; \Omega(t) \, ,
 \label{eq:CG}
\end{equation}
with boundary conditions
\begin{equation}
 u \;=\; 0 \quad \;\;\mbox{and}\;\; \quad \nabla u \cdot \mathbf{n} \;=\; 0 \hspace{10mm}\quad\qquad \mbox{on} \;\;\; \partial \Omega_M(t) \, ,
 \label{eq:cg_bcs}
\end{equation}
along with appropriate Dirichlet or Neumann conditions on any part of the boundary which is fixed.
Note that the moving boundary conditions are equivalent to $u \, \mathbf{v} \cdot \mathbf{n} = - \nabla u \cdot \mathbf{n}$ (a kinematic condition representing zero flux through the moving boundary) and $u = 0$.

\subsubsection{Algorithm for the Crank-Gupta Problem}

 We choose $\mathbb{M}(u) = u$ with ${\cal L}u \equiv \nabla^2 u - 1$
and in the derivation of the discrete forms we define $(\mathcal{L} u)_h$ using
\begin{equation}
 \int_{\Omega_h} w_h \, (\mathcal{L} u)_h \; \mathrm{d}\mathbf{x}
   \;=\; \int_{\partial \Omega_h} w_h \, \nabla u_h \cdot \mathbf{n} \; \mathrm{d}s
     \;-\; \int_{\Omega_h} \nabla w_h \cdot \nabla u_h \; \mathrm{d}\mathbf{x}
     \;-\; \int_{\Omega_h} w_h \; \mathrm{d}\mathbf{x} \, . 
\end{equation}
As in \cite{BHJ05a,BHJ11},
the velocity potential is found using Equations \eqref{eq:potentialequation_discrete} which, when
combined with the boundary conditions, $u = \nabla u \cdot \mathbf{n} = 0$ on
$\partial \Omega_M(t)$
and $\mathbf{v} \cdot \mathbf{n} = 0$ on $\partial \Omega_F$, gives
\begin{eqnarray}
 & & c_{h,u}(w_h) \, \dot{\theta}_{h,u}(t)
     \;+\; \int_{\Omega_h(t)} u_h \, \nabla{\phi}_h \cdot \nabla w_h \; \mathrm{d}\mathbf{x} \nonumber \\
 &=& \int_{\partial \Omega_{h,F}(t)} w_h \, \nabla u_h \cdot \mathbf{n} \; \mathrm{d}s
                   \;-\; \int_{\Omega_h(t)} \nabla w_h \cdot \nabla u_h \; \mathrm{d}\mathbf{x}
                      \;-\; \int_{\Omega_h(t)} w_h \; \mathrm{d}\mathbf{x} \, ,
 \label{eq:cg_phi}
\end{eqnarray}
after integration by parts, along with
\begin{equation}
 \dot{\theta}_{h,u}(t) \;=\;
   \int_{\partial \Omega_{h,F}(t)} \nabla u_h \cdot \mathbf{n} \; \mathrm{d}s
              \;-\; \int_{\Omega_h(t)} \; \mathrm{d}\mathbf{x} \, .
 \label{eq:cg_thetadot}
\end{equation}
Note that boundary integrals only appear on fixed boundaries.
Similarly, the ALE equation (\ref{eq:ALEupdate_discrete}) and an interpolated mesh velocity,
after integration by parts and imposition of the same boundary conditions, gives
\begin{equation}
 \dot{\mu}_{h,u}(t;\hat{w}_h) \;=\; 
    \int_{\partial \Omega_{h,F}} \hat{w}_h \, \nabla u_h \cdot \mathbf{n} \; \mathrm{d}s
      \;-\; \int_{\Omega_h(t)} \nabla \hat{w}_h \cdot
                     ( \nabla u_h + u_h \, \tilde{\mathbf{v}}_h ) \; \mathrm{d}\mathbf{x} \, 
                     \;-\; \int_{\Omega_h(t)} \hat{w}_h \; \mathrm{d}\mathbf{x} \, ,
 \label{eq:CGALE}
\end{equation}
where the $\hat{w}_h$ are the modified test functions of \cite{HBJ09}.

\subsubsection{Numerical Results}

 In \cite{BCR75} an exact solution to the one-dimensional equation on the interval
$[0,X(t)]$ with the boundary conditions in Equation \eqref{eq:cg_bcs} imposed at the moving
boundary is given as
\begin{equation}
 u(x,t) \;=\; \left\{ \begin{array}{ll}
                - x - t + e^{x+t-1} & \;\;\; x \le 1-t \, , \\
                                  0 & \;\;\; x  >  1-t \, .
                      \end{array} \right.
 \label{eq:CGsol}
\end{equation}
The initial condition for our numerical experiments is taken at $t = 0$, when
$u(x,0)= - x + e^{x-1}$ for $x \in [0,1]$. A Neumann condition is imposed at the
fixed boundary at $x = 0$, where $u_x(0,t)= - 1 + e^{t-1}$. A time-step of
$\Delta t = 1.25 \times 10^{-5}$ is used on the coarsest (10-element) mesh.

\begin{figure}[htbp]
 \centerline{\includegraphics[width=0.65\textwidth]{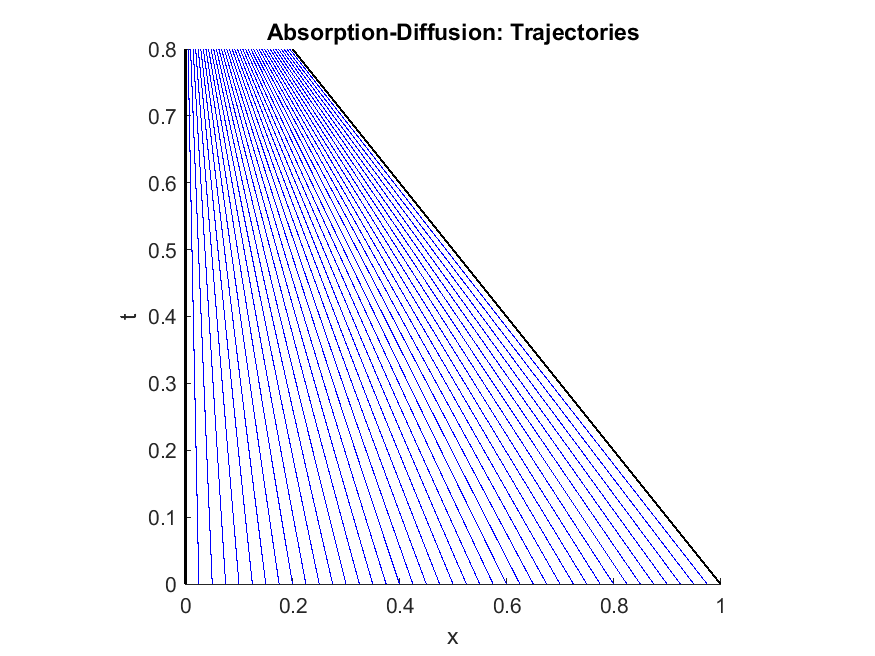} \hspace{-19mm}
          \includegraphics[width=0.65\textwidth]{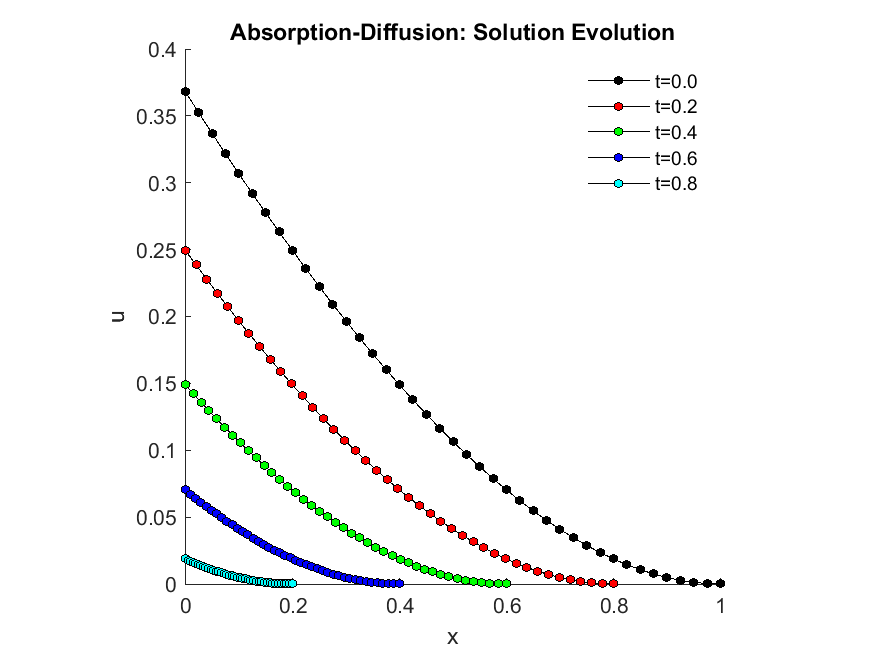}}
 \caption{One-dimensional Crank-Gupta problem approximated on a uniform 40-element mesh:
          node trajectories (left); snapshots of solution (right).}
 \label{fig:CG_evolution}
\end{figure}

\begin{figure}[htbp]
 \centerline{\includegraphics[width=0.65\textwidth]{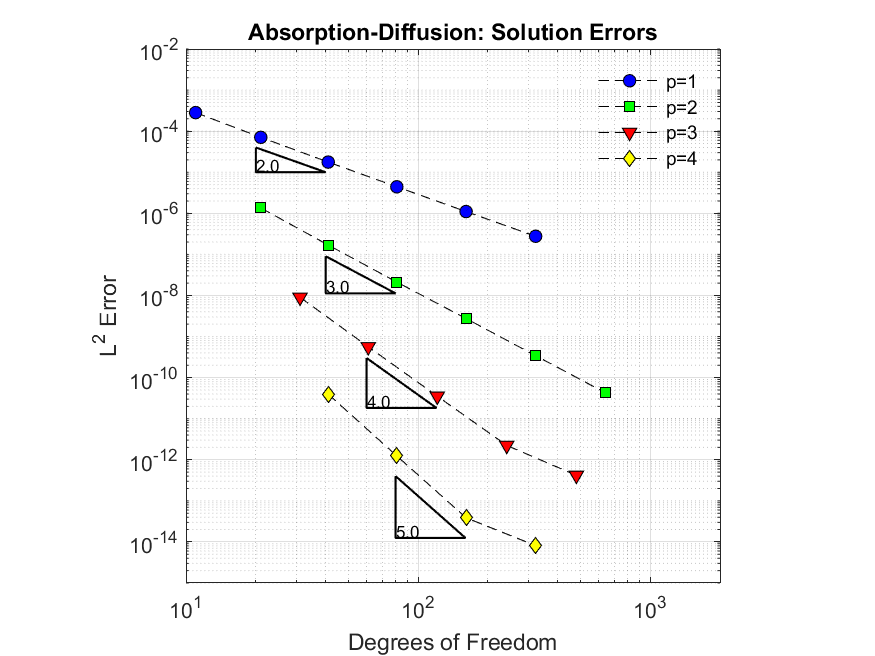} \hspace{-19mm}
          \includegraphics[width=0.65\textwidth]{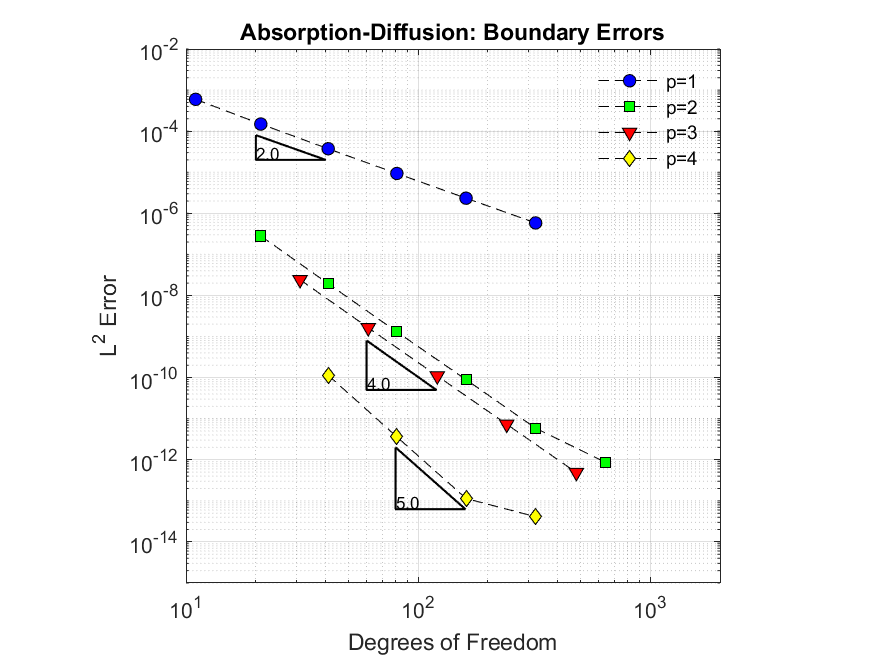}}
 \caption{Mesh convergence plots for the one-dimensional Crank-Gupta problem approximated
          on the interval $t \in [0,0.3]$: solution error (left); boundary position error (right).
          The graphs show results for uniform initial meshes.}
 \label{fig:CG_errors}
\end{figure}

 The evolution of the node positions and the solution can be seen in Figure
\ref{fig:CG_evolution} for the 40-element mesh.
The graphs in Figure \ref{fig:CG_errors} show the expected orders of convergence for $p = 1,2,3,4$
when considering errors in the solution. The same patterns are seen for the errors in the moving
boundary position, except for $p = 2$, which appears to converge at the same rate as $p = 3$.
The reason for this is not certain, but we note three points which may be relevant: 
the exact boundary position in Equation \eqref{eq:CGsol} varies linearly in time;
for $p \geq 2$, both moving boundary conditions in Equation \eqref{eq:cg_bcs}
can be represented exactly by the approximation; the convergence behaviour imitates
that of the time-stepping scheme when $\Delta t \propto h^2$ and the spatial error
is negligible. The same rate of convergence is seen when the internal nodes of the
initial meshes are randomly perturbed.
The effects of rounding error can be seen for higher values of $p$ on finer meshes, beyond which the
error starts to increase with the number of degrees of freedom.

\subsection{Flow in a Porous Medium}
\label{sec:PME}

 The porous medium equation (PME) models flow through a porous medium and provides a nonlinear
PDE for which mass is conserved exactly. It has a kinematic boundary condition but is selected here to
provide a more challenging (nonlinear) test for the moving-mesh framework as applied to an
implicit moving boundary problem where local mass conservation (via the choice of $\mathbb{M}(u) = u$)
is used to determine the boundary motion.

\subsubsection{The Porous Medium Equation}

The porous medium equation, details of which can be found in \cite{Aro86,Vaz06}, is given by
\begin{equation}
 \frac{\partial u}{\partial t} \;=\; \nabla \cdot ( u^m \nabla u )
   \hspace*{10mm} {\mbox{in}} \;\;\; \Omega(t) \, ,
 \label{eq:pme}
\end{equation}
where $m>0$ is usually chosen to be an integer exponent. The boundary conditions are
\begin{equation}
 u \;=\; 0 \quad \;\;\mbox{and}\;\; \quad u \, \mathbf{v} \cdot \mathbf{n} \;=\; - u^m \nabla u \cdot \mathbf{n} \quad\qquad \mbox{on} \;\;\; \partial \Omega_M(t) \, , \\
 \label{eq:pme_bcs}
\end{equation}
along with appropriate Dirichlet or Neumann conditions on any part of the boundary which is fixed.
In the following, we assume that $\partial \Omega_F = \emptyset$, which is the case in all of the
test cases we investigate. A consequence of this and the boundary conditions given in Equation
\eqref{eq:pme_bcs} is that no boundary integrals appear in the following description of the algorithm.

\subsubsection{Algorithm for the Porous Medium Equation}

 We choose $\mathbb{M}(u) = u$ and ${\cal L}u \equiv \nabla \cdot ( u^m \nabla u )$,
as in \cite{BHJ05a,BHJ11},
and in the derivation of the discrete forms we define $(\mathcal{L} u)_h$ using
\begin{equation}
 \int_{\Omega_h} w_h \, (\mathcal{L} u)_h \; \mathrm{d}\mathbf{x}
   \;=\; \int_{\partial \Omega_h} w_h \, (u_h)^m \, \nabla u_h \cdot \mathbf{n} \; \mathrm{d}s
     \;-\; \int_{\Omega_h} (u_h)^m \, \nabla w_h \cdot \nabla u_h \; \mathrm{d}\mathbf{x} \, .
\end{equation}
The velocity potential
is found using Equations \eqref{eq:potentialequation_discrete} which, after integration
by parts and applying the boundary conditions in \eqref{eq:pme_bcs}, give
\begin{equation}
  c_{h,u}(w_h) \, \dot{\theta}_{h,u}(t) \,+\,
           \int_{\Omega_h(t)} u_h \, \nabla \phi_h \cdot \nabla w_h \; \mathrm{d}\mathbf{x}
    \;=\; -\, \int_{\Omega_h(t)} (u_h)^m \, \nabla w_h \cdot \nabla u_h \; \mathrm{d}\mathbf{x} \, ,
 \label{eq:pme_phi}
\end{equation}
along with $\dot{\theta}_{h,u}(t) = 0$, \emph{i.e.}\ mass is conserved. As a result, the first term in Equation
\eqref{eq:pme_phi} can be ignored and it is not necessary to
compute $c_{h,u}(w_h)$ in this case. Similarly, the ALE equation \eqref{eq:ALEupdate_discrete}, after integration by
parts and imposition of the same boundary conditions, gives
\begin{equation}
 \dot{\mu}_{h,u}(t;\hat{w}_h) \;=\; -\, \int_{\Omega_h(t)} \nabla \hat{w}_h \cdot \{ (u_h)^m \, \nabla u_h + u_h \, \tilde{\mathbf{v}}_h \} \; \mathrm{d}\mathbf{x} \, ,
 \label{eq:pmeALE}
\end{equation}
where $\hat{w}_h$ are the modified test functions defined in \cite{HBJ09} and $\tilde{\mathbf{v}}_h$
is a piecewise linear velocity field.

\subsubsection{Numerical Results}

 Equation (\ref{eq:pme}) admits a family of exact compactly supported similarity solutions with moving
boundaries on which $u=0$ \cite{B03,M02}. These were used to confirm the convergence rates of the
second-order method but are of less use for higher orders: (i) when $m = 1$ in \eqref{eq:pme} the similarity
solution is quadratic in space so when $p \geq 2$ the spatial variation can be represented exactly
by the approximation; (ii) when $m > 1$ the similarity solution has reduced regularity at the
moving boundary, so optimal convergence rates are not achieved on uniform initial meshes (a
problem inherited from sub-optimal convergence rates of the approximation of the initial conditions).

\begin{figure}[htbp]
 \centerline{\includegraphics[width=0.65\textwidth]{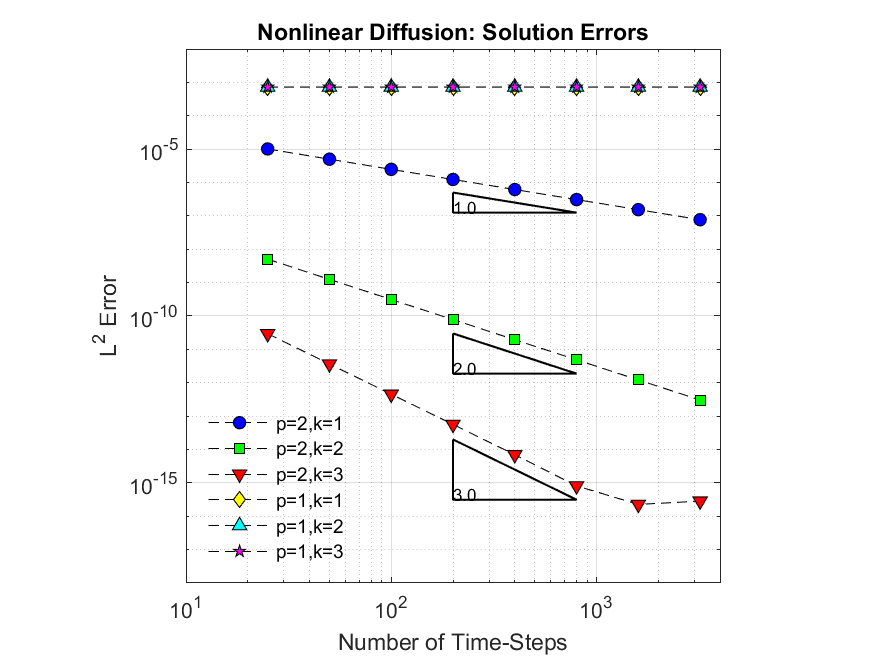} \hspace{-19mm}
          \includegraphics[width=0.65\textwidth]{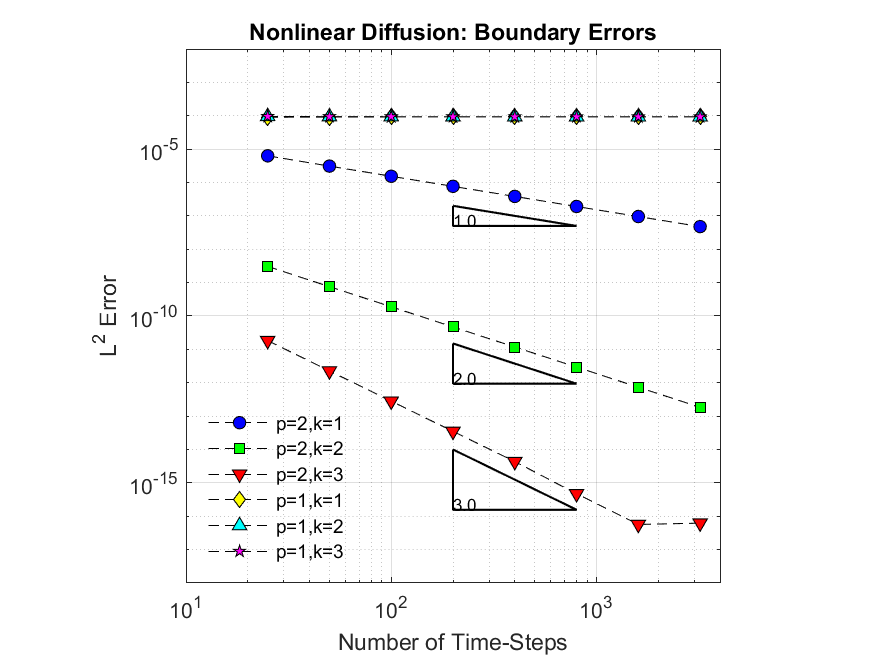}}
 \caption{Mesh convergence plots for the one-dimensional porous medium equation approximated
          on the interval $t \in [t^0,t^0+0.01]$: solution error (left); boundary error (right).
          The graphs show results for uniform initial meshes with spatial polynomial degree $p$
          and temporal approximation order of accuracy $k$.}
 \label{fig:PME_errors_ts}
\end{figure}

 The method is first tested with the one-dimensional similarity solution \cite{B03,M02}:
\begin{equation}
  u(x,t) \;=\; \left\{ \begin{array}{ll}
     \frac{1}{\lambda} \left( 1 - ( \frac{x}{x_0 \lambda} )^{2} \right)^{\frac{1}{m}} \quad & |x| \leq x_0 \lambda \, , \\
                            0 & |x| > x_0 \lambda \, , \end{array} \right.
\end{equation}
in which $\lambda = ( t / t^0 )^{1/(2+m)}$ and $t^0 = 0.5 \, m \, {x_0}^2 / (m+2)$, where we choose $m = 1$ and $x_0 = 0.5$.
Figure \ref{fig:PME_errors_ts} shows the errors obtained when linear ($p=1$) and quadratic ($p=2$)
spatial approximations are used on a mesh of 10 elements with time-stepping methods of three different
orders ($k = 1,2,3$). For $p = 1$ the results are dominated by the spatial approximation error so,
since the same initial spatial mesh is used for all the tests, the total errors remain almost unchanged as the time-step is
decreased. When $p=2$ the spatial variation of the similarity solution can be represented exactly
and the errors, for both the solution and the boundary position, decrease in accordance with the order
of accuracy of the time-stepping method used as the time-step is decreased.
The results are not presented for $p > 2$ because they overlay those for $p = 2$.
The observed convergence rates for $p \geq 2$ match those expected of the
time-stepping scheme alone without the need for any modification to satisfy a
discrete geometric conservation law.

 To confirm the orders of accuracy for $p \geq 2$, we use a problem for which
there is not an available exact solution.
The initial condition is taken to be $u(x,0) = \cos ( \pi x )$
for $x \in [-0.5,0.5]$, so there are moving boundaries with $u = 0$ at each end of the
domain. A time-step of $\Delta t = 1.0 \times 10^{-4}$
is used on the coarsest (10-element) mesh.

 Since the exact solution is not available, we investigate convergence rates by comparing approximations at
successive levels of refinement, \emph{i.e.}\ instead of computing the error norms in Equation \eqref{eq:analyticerror}
we consider mixed discrete/continuous norms
\begin{eqnarray}
  \mbox{Error}_u &\approx& \left( \sum_{m=1}^{N_t'} \Delta t' \int_{\Omega_h(t^m)} [ u_{2h}(\mathbf{x},t^m) - u_h(\mathbf{x},t^m) ]^2 \; \mathrm{d}\mathbf{x} \right)^{\frac{1}{2}} \, , \nonumber \\
  \mbox{Error}_x &\approx& \left( \sum_{m=1}^{N_t'} \Delta t' \, | \mathbf{x}_{\partial,2h}(t^m) - \mathbf{x}_{\partial,h}(t^m) |^2 \right)^{\frac{1}{2}} \, , 
  \label{eq:approximateerror}
\end{eqnarray}
in which $N_t'$ is a specified number of sample points in time, each $\Delta t'$ apart, and the subscript
$\cdot_{2h}$ indicates a quantity computed for a simulation in which the initial mesh size is double that
of a simulation denoted by the subscript $\cdot_{h}$.
This difference should converge at the same rate as the exact error, but its convergence to zero
does not guarantee convergence to the correct solution.
 In the following simulations, $\Delta t' = 1.0 \times 10^{-4}$ is used in the error estimates
(to match the time-step used on the coarsest mesh
and give discrete time points at which approximations are computed in all of the simulations).
This provides $N_t' =100$ uniformly distributed time points at which the error is sampled in the simulation interval, $t \in [0,0.01]$.

\begin{figure}[htbp]
 \centerline{\includegraphics[width=0.65\textwidth]{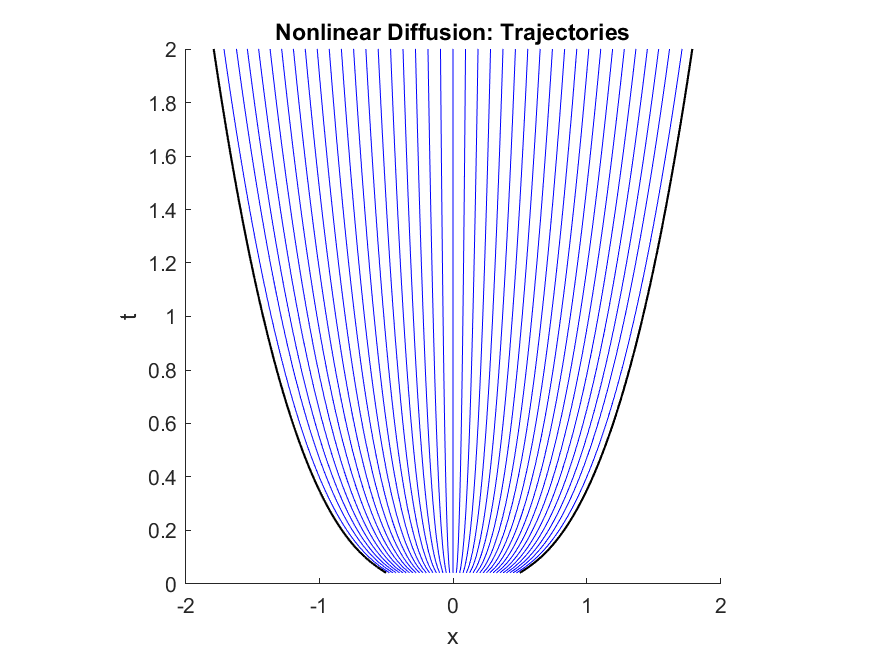} \hspace{-19mm}
          \includegraphics[width=0.65\textwidth]{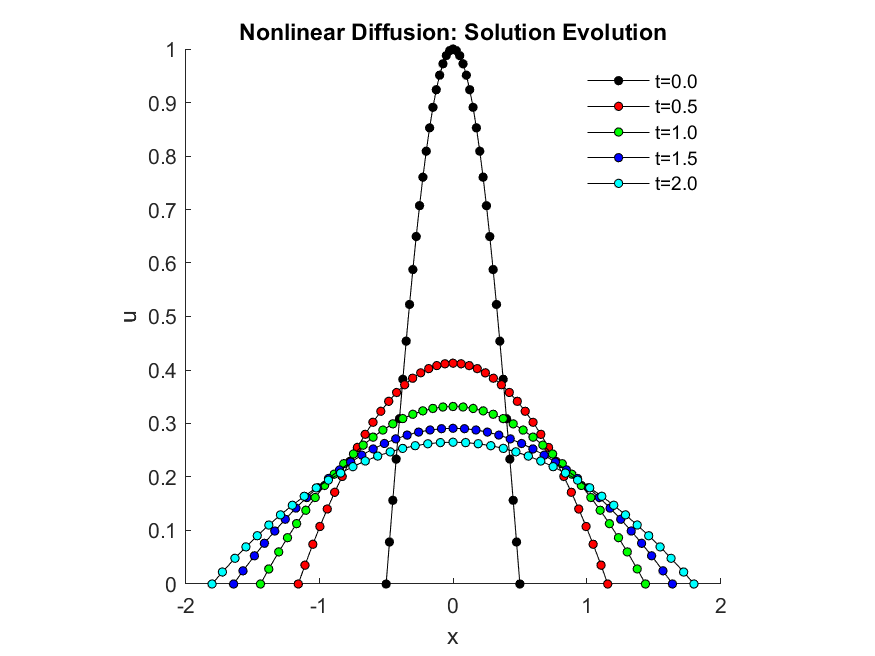}}
 \caption{One-dimensional porous medium equation approximated on a uniform 40-element mesh:
          node trajectories (left); snapshots of solution (right).}
 \label{fig:PME_evolution}
\end{figure}

\begin{figure}[htbp]
 \centerline{\includegraphics[width=0.65\textwidth]{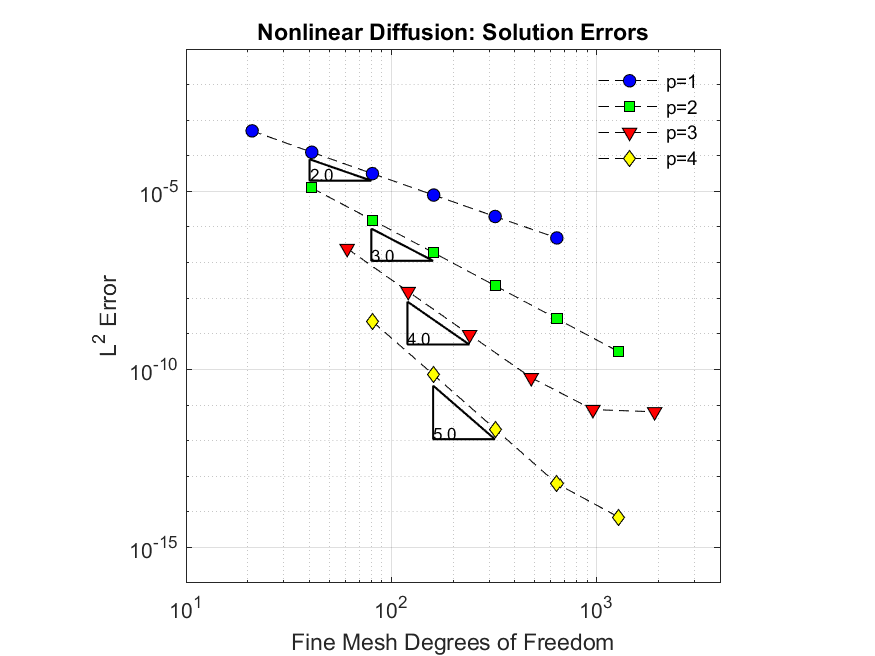} \hspace{-19mm}
          \includegraphics[width=0.65\textwidth]{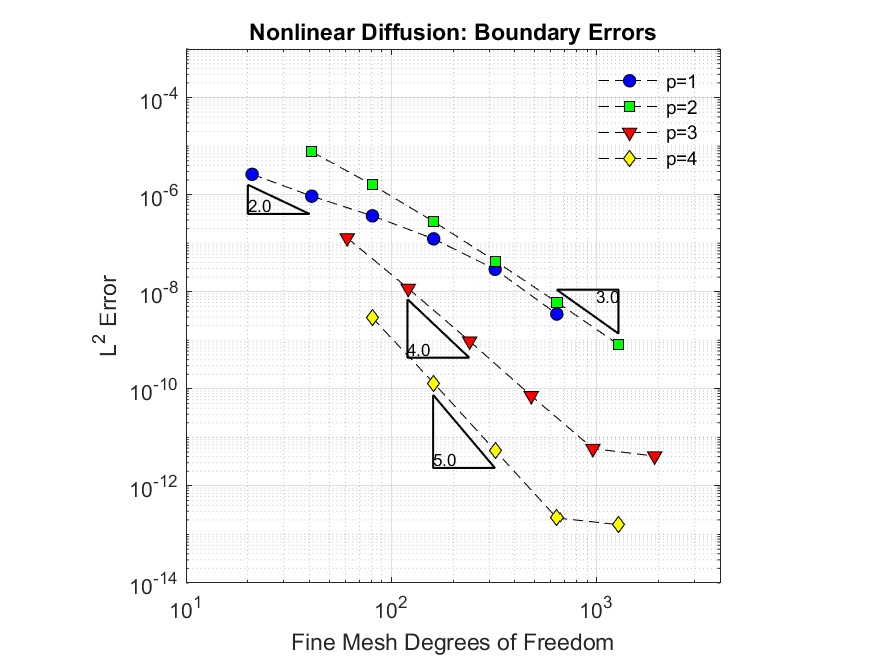}}
 \caption{Mesh convergence plots for the one-dimensional porous medium equation approximated
          on the interval $t \in [0,0.01]$: solution error (left); boundary error (right).
          The graphs show results for uniform initial meshes.}
 \label{fig:PME_errors}
\end{figure}

 The evolution of the node positions and the solution can be seen in Figure \ref{fig:PME_evolution}
for the 40-element mesh.
The graphs in Figure \ref{fig:PME_errors} show the expected orders of convergence for $p = 1,2,3,4$
when considering errors in the solution. The same patterns are seen for the errors in the moving
boundary position, except for $p = 1$, which appears to converge at the same rate as $p = 2$.
The reason for this is not clear but it may be an artefact of the uniform initial mesh since this apparent
superconvergence disappears when the interior nodes of the initial mesh are perturbed.
The effects of rounding error can be seen for higher values of $p$ on finer meshes, beyond which the
error starts to increase with the number of degrees of freedom.

\section{Discussion and Conclusions}
\label{sec:conc}

 In this paper a moving-mesh finite element method has been described which
achieves high-order accuracy in approximating the evolution of both the
solution and the boundary position for a range of one-dimensional parabolic implicit moving boundary
problems. Fifth-order accuracy has been demonstrated for linear diffusion with
phase change across a moving interface, linear diffusion with absorption, and nonlinear
diffusion. The method is explicit in its prediction of both the mesh positions and solution
values at each new time level (so only linear systems need to be inverted within each
time-step, even for nonlinear problems) and requires no mapping to a reference domain
or remapping of the solution from one mesh to the next. This is the first time that such
high orders of accuracy have been achieved for implicit moving boundary problems within
a conservative ALE framework and it has been done without the need to modify the time-stepping
schemes to satisfy a discrete geometric conservation law.

 The key observation that allows the moving-mesh algorithm to inherit the order of
accuracy of the underlying fixed-mesh method is that high-order approximation of mesh
velocities is only required at the boundary.
 The representation of the velocity field
in the interior of the domain does not have to be high-order: it can and, where possible, should be piecewise linear
(as also suggested in \cite{LXY23} for prescribed boundary motion in multiple space dimensions).
As a consequence, the assumption that the finite element test functions are transported with the
mesh velocity field (which is required for this conservative, mapping-free approach)
means that the standard finite element basis functions on a mesh are transported to
the corresponding standard finite element basis functions on the mesh at the next time level, so
no remapping or projection step is required.

 For moving boundary problems, the selection of a piecewise linear velocity field
does not generalise straightforwardly to multiple space dimensions, due to the presence
of curved boundary faces. If the boundary is fixed, and can be represented exactly with
flat faces, then meshes of simplices in the interior can still be moved with a piecewise
linear velocity field which transports the basis functions in an appropriate manner.
We anticipate that the multidimensional framework presented in Section
\ref{sec:algorithm} would achieve high-order accuracy for fixed-boundary problems where
an internal mesh consisting of simplicial elements is
moved. However, this is not our immediate focus because our goal is to develop a mapping-free
approach in the physical domain which retains high-order accuracy for general (curved)
multidimensional moving boundary problems.

 From a practical perspective, the algorithm would benefit from
modifications which could make it more robust and efficient, particularly for supporting future
developments in multiple space dimensions. For the parabolic problems considered
here, the stability constraints inherent to explicit time-stepping algorithms can make the time-step
prohibitively small. It would be simple to switch to an implicit time-stepping scheme in Section
\ref{sec:algorithm} if the mesh motion were prespecified. However, it is not clear how to retain
the benefits of unconditionally stable time-stepping when the positions of the mesh nodes
and the domain boundary at the new time level depend on the solution.
Some progress on this has been made for a finite difference approach in one dimension \cite{BL14}.

 In this work we have consciously not investigated alternative approaches to determining the
internal mesh node velocities. When the mass monitor is used, it generates accurate predictions
of the boundary movement, but it creates a Lagrangian method which inherits
the well-known robustness issues which arise when the mesh follows the flow. However, the
internal mesh movement could be generated by any method, {\emph e.g.}\ one which adapts to a
local error indicator and includes regularisation to avoid mesh tangling and optimise mesh quality.
As long as the chosen technique is used to generate piecewise linear mesh velocity fields,
the conservative ALE step in Section \ref{sec:algorithm} should retain the order
of accuracy of the fixed mesh method used to approximate the PDE, without the need for
higher-order quadrature.

 Finally, the robustness and accuracy of mesh movement algorithms can often be improved
by changing the mesh structure. This might involve edge/face swapping or insertion, to
locally improve mesh quality, or $h$-refinement, guided by a local error indicator or
estimate. Theoretical analysis might also help to indicate when mesh movement would be
more effective than $h$- or $p$-adaptivity. If robust methods can be developed which retain
high-order accuracy on moving meshes then it opens up the possibility of $hpr$-adaptive
algorithms for computationally efficient simulation of implicit moving boundary problems.

\section*{Acknowledgements}

 T.\ J.\ Radley is grateful for the financial support of both EPSRC and of the European Union (ERC Synergy, NEMESIS, project number 101115663).
Views and opinions expressed are however those of the authors only and do not necessarily reflect those of the European Union or the European
Research Council Executive Agency. Neither the European Union nor the granting authority can be held responsible for them.

\end{document}